\documentclass{ws-idaqp}

\usepackage[latin1]{inputenc}

\usepackage{amsmath}
\usepackage{amssymb}
\usepackage[mathscr]{eucal}

    \usepackage{theorem}

\newenvironment{env}[2]{\begin{#1}#2\end{#1}}{}
    \newcommand{\beq}[1]{\begin{env}{equation}{#1}}
    \newcommand{\beqn}[1]{\begin{env}{equation*}{#1}}
    \newcommand{\bal}[1]{\begin{env}{align}{#1}}
    \newcommand{\baln}[1]{\begin{env}{align*}{#1}}
    \newcommand{\bga}[1]{\begin{env}{gather}{#1}}
    \newcommand{\bgan}[1]{\begin{env}{gather*}{#1}}
    \newcommand{\bflal}[1]{\begin{env}{flalign}{#1}}
    \newcommand{\bflaln}[1]{\begin{env}{flalign*}{#1}}
    \newcommand{\bmu}[1]{\begin{env}{multline}{#1}}
    \newcommand{\bmun}[1]{\begin{env}{multline*}{#1}}
    \newcommand{\bsp}[1]{\begin{env}{split}{#1}}

    \newcommand{\eeq}{\end{env}}
    \newcommand{\eeqn}{\end{env}}
    \newcommand{\eal}{\end{env}}
    \newcommand{\ealn}{\end{env}}
    \newcommand{\ega}{\end{env}}
    \newcommand{\egan}{\end{env}}
    \newcommand{\eflal}{\end{env}}
    \newcommand{\eflaln}{\end{env}}
    \newcommand{\emu}{\end{env}}
    \newcommand{\emun}{\end{env}}
    \newcommand{\esp}{\end{env}}

\newcommand{\lf}{\vspace{2ex}}

\renewcommand{\bf}[1]{\textbf{#1}}
\renewcommand{\it}[1]{\textit{#1}}
\renewcommand{\sc}[1]{\textsc{#1}}
\renewcommand{\sf}[1]{\textsf{#1}}
\renewcommand{\sl}[1]{\textsl{#1}}
\renewcommand{\tt}[1]{\texttt{#1}}
\newcommand{\hl}[1]{\bf{\it{#1}}}

\newcommand{\mbf}[1]{\mathbf{#1}}
\newcommand{\msf}[1]{\text{\small$\sf{#1}$}}

\newcommand{\cmc}[1]{\mathcal{#1}}
\newcommand{\eus}[1]{\mathscr{#1}}
\newcommand{\euf}[1]{\mathfrak{#1}}
\newcommand{\bb}[1]{\mathbb{#1}}

\newcommand{\nbd}[1]{$#1$\nobreakdash--}
\newcommand{\ol}[1]{\overline{#1}}
\newcommand{\ul}[1]{\underline{#1}}

\newcommand{\vt}{\vartheta}

\newcommand{\norm}[1]{\left\lVert#1\right\rVert}

\newcommand{\bfam}[1]{\bigl(#1\bigr)}

\newcommand{\AB}[1]{\langle#1\rangle}

\newcommand{\CB}[1]{\{#1\}}
\newcommand{\bCB}[1]{\bigl\{#1\bigr\}}
\newcommand{\BCB}[1]{\Bigl\{#1\Bigr\}}

\newcommand{\sbars}[1]{\:\bar{#1}^s\:}
\newcommand{\sodots}{\sbars{\odot}}
\newcommand{\set}[2][]{
    \ifthenelse{\equal{#1}{}}{
        \CB{#2}}{
        \CB{#1~|~#2}}}
\newcommand{\bset}[2][]{
    \ifthenelse{\equal{#1}{}}{
        \bCB{#2}}{
        \bCB{#1~|~#2}}}
\newcommand{\Bset}[2][]{
    \ifthenelse{\equal{#1}{}}{
        \BCB{#2}}{
        \BCB{#1~\big|~#2}}}

\DeclareMathOperator{\ls}{\normalfont\msf{span}}
\DeclareMathOperator{\cls}{\ol{\ls}}

\DeclareMathOperator{\id}{\normalfont\msf{id}}

\newcommand{\C}{\bb{C}}

\newcommand{\E}{\bb{E}}

\newcommand{\N}{\bb{N}}

\newcommand{\R}{\bb{R}}

\newcommand{\cA}{\cmc{A}}
\newcommand{\cB}{\cmc{B}}
\newcommand{\cC}{\cmc{C}}

\newcommand{\cK}{\cmc{K}}

\newcommand{\sB}{\eus{B}}

\newcommand{\sF}{\eus{F}}

\newcommand{\sK}{\eus{K}}

\newcommand{\el}{\euf{l}}

\newcommand{\er}{\euf{r}}

\newcommand{\eC}{\euf{C}}

\newcommand{\eW}{\euf{W}}

\newcommand{\U}{\mbf{1}}

\newcommand{\f}{\text{\scriptsize$\sF$}}

    \numberwithin{equation}{section}
    \renewcommand{\theequation}{\thesection.\arabic{equation}}

        \newcommand{\mnname}{Mathematical note.}
        \newcommand{\enname}{End of the note.}
        \newcommand{\definame}{Definition.}
        \newcommand{\propname}{Proposition.}
        \newcommand{\lemname}{Lemma.}
        \newcommand{\exname}{Example.}
        \newcommand{\exername}{Exercise.}
        \newcommand{\remname}{Remark.}
        \newcommand{\obname}{Observation.}
        \newcommand{\thmname}{Theorem.}
        \newcommand{\corname}{Corollary.}

        \renewcommand{\proof}{\lf\noindent\bf{Proof.~}}

    \theoremheaderfont{\normalfont\bfseries}
    \theoremstyle{change}
        \theorembodyfont{\rmfamily}
            \newtheorem{emp}{}[section]
                \newcommand{\bemp}[1][]{
                    \begin{emp}\hskip-\labelsep\bf{#1}\hskip\labelsep}
                \newcommand{\eemp}{\end{emp}}
\newtheorem{itemp}[emp]{}
                \newcommand{\bitemp}[1][]{
                    \begin{itemp}\hskip-\labelsep\bf{#1}\hskip\labelsep\normalfont\itshape}
                \newcommand{\eitemp}{\end{itemp}}
            \newtheorem{mn}[emp]{\mnname}
                \newcommand{\bmn}{\begin{mn}~\begin{quotation}\renewcommand{\baselinestretch}{1}\small\noindent\ignorespaces}
                \newcommand{\emn}{\end{quotation}\hfill\bf{\enname}\end{mn}}
            \newtheorem{ex}[emp]{\exname}
                \newcommand{\bex}{\begin{ex}}
                \newcommand{\eex}{\end{ex}}
            \newtheorem{exer}[emp]{\exername}
                \newcommand{\bexer}{\begin{exer}}
                \newcommand{\eexer}{\end{exer}}
            \newtheorem{defi}[emp]{\definame}
                \newcommand{\bdefi}{\begin{defi}}
                \newcommand{\edefi}{\end{defi}}
            \newtheorem{rem}[emp]{\remname}
                \newcommand{\brem}{\begin{rem}}
                \newcommand{\erem}{\end{rem}}
            \newtheorem{ob}[emp]{\obname}
                \newcommand{\bob}{\begin{ob}}
                \newcommand{\eob}{\end{ob}}
        \theorembodyfont{\normalfont\itshape}
            \newtheorem{thm}[emp]{\thmname}
                \newcommand{\bthm}{\begin{thm}}
                \newcommand{\ethm}{\end{thm}}
            \newtheorem{prop}[emp]{\propname}
                \newcommand{\bprop}{\begin{prop}}
                \newcommand{\eprop}{\end{prop}}
            \newtheorem{cor}[emp]{\corname}
                \newcommand{\bcor}{\begin{cor}}
                \newcommand{\ecor}{\end{cor}}
            \newtheorem{lem}[emp]{\lemname}
                \newcommand{\blem}{\begin{lem}}
                \newcommand{\elem}{\end{lem}}
\newenvironment{empn}[1]{\lf\noindent\bf{#1}\ignorespaces\hskip\labelsep}{\lf}
		\newcommand{\bempn}[1]{\begin{empn}{#1}}
		\newcommand{\eempn}{\end{empn}}
		\newcommand{\bitempn}[1]{\begin{empn}{#1}\normalfont\itshape}
		\newcommand{\eitempn}{\end{empn}}
                \newcommand{\bmnn}{\begin{empn}{\mnname}~\begin{quotation}\renewcommand{\baselinestretch}{1}\small\noindent\ignorespaces}
                \newcommand{\emnn}{\end{quotation}\hfill\bf{\enname}\end{empn}}
		\newcommand{\bexn}{\begin{empn}{\exname}}
		\newcommand{\eexn}{\end{empn}}
		\newcommand{\bexern}{\begin{empn}{\exername}}
		\newcommand{\eexern}{\end{empn}}
		\newcommand{\bdefin}{\begin{empn}{\definame}}
		\newcommand{\edefin}{\end{empn}}
		\newcommand{\bremn}{\begin{empn}{\remname}}
		\newcommand{\eremn}{\end{empn}}
		\newcommand{\bobn}{\begin{empn}{\obname}}
		\newcommand{\eobn}{\end{empn}}

\renewcommand{\msf}[1]{\mathsf{#1}}

\renewcommand{\thefootnote}{[\alph{footnote}]}

\begin{document}

\title{Representations of $\sB^a(E)$\renewcommand{\thefootnote}{}\thanks{2000 AMS-Subjecrt classification: 46L08, 46L55, 46L53, 46M05, 46M15}}

\author{Paul S.\ Muhly\thanks{PSM is supported by grants from the National Science Foundation.}}

\address{Department of Mathematics, University of Iowa\\Iowa City, IA 52242, USA\\E-mail: \tt{pmuhly@math.uiowa.edu}}

\author{Michael Skeide\thanks{MS is supported by DAAD and by research funds of the Department S.E.G.e S.\ of University of Molise.}}

\address{Dipartimento S.E.G.e S., Universit\`a\ degli Studi del Molise, Via de Sanctis\\86100 Campobasso, Italy\\E-mail: \tt{skeide@math.tu-cottbus.de}}

\author{Baruch Solel\thanks{BS is supported by the Fund for the Promotion of Research at the Technion.}}

\address{Department of Mathematics, Technion\\Haifa 32000, Israel\\E-mail: \tt{mabaruch@techunix.technion.ac.il}}

\maketitle


\begin{abstract}
\noindent
We present a short and elegant proof of the complete theory of strict representations of the algebra $\sB^a(E)$ of all adjointable operators on a Hilbert \nbd{\cB}module $E$ by operators on a Hilbert \nbd{\cC}module $F$. An analogue for \nbd{W^*}modules and normal representations is also proved. As an application we furnish a new proof of Blecher's \it{Eilenberg-Watts theorem}.
\end{abstract}





\lf\lf\noindent
It is well known that if $G$ is a Hilbert space and if $\vt\colon\sB(G)\rightarrow\sB(H)$ is a normal unital (\nbd{*})representation of $\sB(G)$ on another Hilbert space $H$, then $H$ factors as $G\otimes H_\vt$ where $H_\vt$ is a third Hilbert space, the \it{multiplicity space}, and $\vt(a)$ acts on $G\otimes H_\vt$ as $a\otimes\id_{H_\vt}$. We have two objectives in this note.  First, we want to show that two theorems of Rieffel's famous paper [\refcite{Rie74}, Theorem 6.23 and Theorem 6.29], suggest a generalization of this result to the setting of representations of $\sB^a(E)$ in $\sB^a(F)$, where $E$ and $F$ are (right) Hilbert (\nbd{C^*})modules over \nbd{C^*}algebras $\cB$ and $\cC$, say, and where $\sB^a(E)$ and $\sB^a(F)$ denote the algebras of adjointable module maps on $E$ and $F$, respectively.  And second, we want to apply this generalization to give a new proof of Blecher's \it{Eilenberg-Watts theorem}[\refcite{Ble97}].  Along the way, we shall see how the generalization fits with familiar facts from Hilbert module theory.

More specifically, for the first objective, we show in Theorem \ref{strirepthm}, which we call \hl{the representation theorem}, that if $\vt$ is a strict unital representation of $\sB^a(E)$ by adjointable operators on a Hilbert \nbd{\cC}module $F$, then $F$ factors into a tensor product $E\odot F_\vt$, where $F_\vt$ is a Hilbert \nbd{\cB}\nbd{\cC}bimodule\footnote{There are several terms in use for bimodules over \nbd{C^*} and \nbd{W^*}algebras. In this note we will use the term \it{correspondence}.}, and $\vt(a)$ acts as $a\odot\id_{F_\vt}$.  Further, in contrast to the Hilbert space setting, where $H_\vt$ is usually viewed as a ``lifeless'' space whose sole purpose is to index the multiplicity of the identity representation in the given representation $\vt$, the space $F_\vt$ has additional internal structure worthy of investigation. In fact, the proof shows that even in the Hilbert space setting,  $H_\vt$ exhibits additional structure. (See Example \ref{B(H)ex}.) With minor adjustments, the representation theorem and the proof apply in the setting of \nbd{W^*}modules over \nbd{W^*}algebras. This \nbd{W^*}representation theorem will be presented in Theorem \ref{normrepthm}.

The representation theorem may be known to the cognoscenti, but it seems to have escaped notice in the literature, although at times, special cases have been proved in \it{ad hoc} ways.  Consequently, this note serves something of a didactic purpose.  Indeed, to make the ideas widely accessible, we review salient aspects of Hilbert modules that we use at the beginning of the next section.\footnote{The reader familiar with Hilbert module theory can skip directly to Remark \ref{Motivation} or to Theorem \ref{strirepthm}, if all that is desired is the statement and (short) proof.}  The surprisingly simple proof rests on the observation (indeed, almost a definition) that if $E$ is a full Hilbert \nbd{\cB}module, then $E$ is a Morita equivalence from $\sK(E)$ to $\cB$.  Under this hypothesis, the \it{multiplicity module} $F_\vt$ of a representation $\vt$ of $\sB^a(E)$ is unique (Theorem \ref{unithm}). In Observation \ref{nonfob} we analyze the extent to which uniqueness fails when $E$ is not necessarily full.  It turns out that a clear understanding of this is essential for the second section, where we apply Theorem \ref{strirepthm} to give a new proof of Blecher's \it{Eilenberg-Watts theorem} [\refcite{Ble97}], which asserts that every (appropriately regular) functor between categories of Hilbert modules is implemented by tensoring with a fixed bimodule. (See Theorem \ref{EWthm}.)  The analysis of what happens when $E$ fails to be full also has relevance for the study of isomorphisms of $\sB^a(E)$ and is related to several notions of Morita equivalence for right modules and bimodules. We will study these ramifications and their connections with the papers of Muhly and Solel [\refcite{MuSo00,MuSo05}] and Skeide [\refcite{Ske03p1}] elsewhere.

The constructions for composed representations iterate associatively (Theorem \ref{assthm}). This has particular relevance for the study of endomorphism semigroups and product systems. We started studying these ramifications and their connections with the papers of Muhly and Solel [\refcite{MuSo00}] and Skeide [\refcite{Ske03p1}] in Skeide [\refcite{Ske04p}].

Our proof of Blecher's Eilenberg-Watts theorem, Theorem \ref{EWthm}, is based entirely upon Theorem \ref{strirepthm} and avoids operator space technology.  In fact, a direct corollary of Theorem \ref{strirepthm} asserts that every representation of $\sB^a(E)$ arises as the restriction of a functor to the object $E$ and its endomorphisms $\sB^a(E)$. We point out in Remark \ref{nbackrem}, however, that Theorem \ref{strirepthm} does not appear to be a consequence of Theorem \ref{EWthm}. In order to avoid the introduction of categories at an early stage, we postpone their use until we really need them --- in the second section.

\section{Representations}\label{RepSEC}

Throughout, $\cA,\cB,\cC,\ldots$ will denote \nbd{C^*}algebras. A \hl{correspondence from $\cB$ to $\cC$} is a (right) Hilbert module $F$ over $\cC$ that is endowed with a left action of $\cB$ via a  nondegenerate representation of $\cB$ in the algebra of adjointable operators on $F$, $\sB^a(F)$.  Rieffel initially introduced the notion of correspondence under the name \nbd{\cC}rigged \nbd{\cB}module [\refcite{Rie74}].  Subsequently, the term Hilbert \nbd{\cB}\nbd{\cC}bimodule has also been used.\footnote{Rieffel's rigged modules are, by definition, full.  It is important for parts of our discussion to allow non-full modules and so, in particular, we do not assume our correspondences are full.}  Every Hilbert \nbd{\cB}module $E$ may be viewed as a correspondence from $\sB^a(E)$ to $\cB$ where $\sB^a(E)$ is the algebra of adjointable operators on $E$. Every \nbd{C^*}algebra $\cB$ may be viewed as a correspondence, the \hl{trivial} correspondence, from $\cB$ to $\cB$ when equipped with the natural bimodule structure coming from right and left multiplication and inner product $\AB{b,b'}=b^*b'$.

\brem
The assumption of nondegeneracy of the left action of a correspondence is crucial in what follows. It eliminates a couple of cases that, otherwise, would cause problems. (Some of these problems can be circumvented, but only if the submodule $\cls\cB F$ of $F$ is complemented in $F$.)
\erem

The algebra $\sK(E)\subset\sB^a(E)$ of \hl{compact operators} is the norm completion  of the algebra $\sF(E)$ of \hl{finite-rank operators}, which is spanned by the \hl{rank-one operators} $xy^*\colon z\mapsto x\AB{y,z}$, $(x,y\in E$). The \hl{dual module} of $E$ is $E^*$, equipped with the \nbd{\sK(E)}valued inner product $\AB{x^*,y^*}=xy^*$ and with its natural \nbd{\cB}\nbd{\sB^a(E)}module structure $bx^*a=(a^*xb^*)^*$, is a correspondence from $\cB$ to $\sB^a(E)$. (Note that the left action of $\cB$ on $E^*$ is nondegenerate because the right action of $\cB$ on $E$ is nondegenerate and the canonical mapping $E\rightarrow E^*$ is an anti-linear isometry of Banach spaces.) Whenever, $\cA\subset\sB^a(E)$ is a \nbd{C^*}algebra which contains $\sK(E)$, then we may view $E^*$ also as a correspondence from $\cB$ to $\cA$. In particular, $E^*$ is also a correspondence from $\cB$ to $\sK(E)$. It follows that the action of $\sK(E)$ on $E$ is nondegenerate (because its action on $E^*$ is), so that when $\sK(E) \subseteq \cA \subseteq \sB^a(E)$, $E$ can also be viewed either as a correspondence from $\sK(E)$ to $\cB$ or from $\cA$ to $\cB$.

$\sB^a(E)$ is the multiplier algebra of $\sK(E)$; see [\refcite{Kas80,Ske01}]. Therefore, $\sB^a(E)$ inherits a \hl{strict} topology. The strict topology coincides with the \nbd{*}strong topology, when restricted to bounded subsets; see [\refcite{Lan95,Ske01}]. We say that a bounded linear mapping $\sB^a(E)\rightarrow\sB^a(F)$ is \hl{strict}, if it is strictly continuous on bounded subsets. Every nondegenerate representation $\sK(E)\rightarrow\sB^a(F)$ extends to a unique strict unital representation $\sB^a(E)\rightarrow\sB^a(F)$; see [\refcite{Lan95}] and also Corollary \ref{NonUnital}. So the assumption that a unital representation $\sB^a(E)\rightarrow\sB^a(F)$ is strict means that it can be reconstructed from its restriction to $\sF(E)$. Clearly, in this case $\sF(E)$ acts nondegenerately on $F$ (because a bounded approximate unit for $\sK(E)$ chosen from the dense subset $\sF(E)$ converges strictly to $\U$). Therefore, a correspondence $F$ from $\sB^a(E)$ to $\cC$ with strict left action may also be viewed as correspondence from $\sK(E)$ to $\cC$ and, conversely, every correspondence from $\sK(E)$ to $\cC$ can be viewed as a correspondence from $\sB^a(E)$ to $\cC$ with strict left action extending that of $\sK(E)$ in a unique way.

\brem
By an application of the \it{closed graph theorem} a linear mapping that is strictly continuous on bounded subsets of $\sB^a(E)$ is bounded.  But usually before we can show that a concrete mapping is strict, we first show that it is bounded and then we check strict continuity on bounded subsets by checking \nbd{*}strong continuity only with respect to a total subset of $F$ (which is not sufficient, if we do not know that the mapping is bounded). Including boundedness into the definition has the advantage that it works nicely also for pre-Hilbert modules. For the same reason, we defined $\sK(E)$ as the completion of $\sF(E)$ instead of its closure in $\sB^a(E)$. ($\sK(E)$ is then always defined as a \nbd{C^*}algebra, even if $E$ is only a pre-Hilbert module. Of course, in this case $\sK(E)$ need not be representable as an algebra of operators on $E$.)
\erem

A correspondence $M$ from $\cB$ to $\cC$ is a \hl{Morita equivalence} (from $\cB$ to $\cC$), if $M$ is \hl{full} (i.e.\ if the range of the inner product of $M$ generates $\cC$), and if the canonical mapping $\cB\rightarrow\sB^a(M)$ is an isomorphism onto $\sK(M)\subset\sB^a(M)$. In the literature, what we are calling a Morita equivalence is also called a \it{strong Morita equivalence}, a \it{Morita-Rieffel equivalence} and an \it{imprimitivity bimodule}.\footnote{Also, in the literature, the notion of Morita equivalence from $\cB$ to $\cC$ is sometimes formulated to include the assertion that the correspondence $M$ has a \nbd{\cB}valued inner product satisfying a compatibility condition with the \nbd{\cC}valued inner product.  However, this condition is covered by the assertion that $\cB$ is isomorphic to $\sK(M)$ [\refcite{Lan95}].  The left \nbd{\cB}valued inner product does not play an explicit role here.} Note in particular that if $E$ is a full Hilbert \nbd{\cB}module, then $E$ is a Morita equivalence from $\sK(E)$ to $\cB$, and $E^*$ is a Morita equivalence from $\cB$ to $\sK(E)$. Even if $E$ is not full, $E$ is still a Morita equivalence from $\sK(E)$ to $\cB_E$, and $E^*$ is a Morita equivalence from $\cB_E$ to $\sK(E)$, where $\cB_E:=\cls\AB{E,E}$ denotes the \hl{range ideal} of $E$ in $\cB$.

The \hl{tensor product} (over $\cB$) of a correspondence $E$ from $\cA$ to $\cB$ and a correspondence $F$ from $\cB$ to $\cC$ is the unique correspondence $E\odot F$ from $\cA$ to $\cC$ that is generated by elementary tensors $x\odot y\in E\odot F$ with inner product
\beqn{
\AB{x\odot y,x'\odot y'}
~=~
\AB{y,\AB{x,x'}y'}.
}\eeqn
In particular, $E\odot E^*$ is identified with $\sK(E)$ (via the map $x\odot y^* \mapsto xy^*$) when viewed as the correspondence from $\sK(E)$ to $\sK(E)$, and $E^*\odot E$ is identified with $\cB_E$ (via the map $x^*\odot y \mapsto \AB{x,y}$).  Consequently, if $M$ is a Morita equivalence from $\cB$ to $\cC$, then under all the identifications that have been made, we may write $M\odot M^*=\cK(M)=\cB$ and $M^*\odot M=\cC$. If $F$ is a correspondence from $\cB$ to $\cC$, then the correspondence $\cC$ serves as right identity under tensor product, i.e.\ $F\odot\cC=F$ (via $y\odot c\mapsto yc$), and (by nondegeneracy) $\cB$ serves as left identity, i.e.\ $\cB\odot F=F$ (via $b\odot y\mapsto by$).

Summarizing: When $F$ is a Hilbert \nbd{\cB}module, then we may construct the Hilbert $\sB^a(E)$--mod\-ule $F\odot E^*$ and, if $E$ is full, we get back $F$ as $F=(F\odot E^*)\odot E$. When $F$ is a Hilbert \nbd{\sK(E)}module, then we may construct the Hilbert \nbd{\cB}module $F\odot E$ and get back $F$ as $F=(F\odot E)\odot E^*$. (Here $E$ need not be full.) Similarly, when $F$ is a correspondence from $\cB$ to $\cC$, then we may construct the correspondence $E\odot F$ from $\sB^a(E)$ to $\cC$ (which has strict left action of $\sB^a(E)$) and, if $E$ is full, we get back $F$ as $F=E^*\odot(E\odot F)$. When $F$ is a correspondence from $\sB^a(E)$ to $\cC$ with strict left action, then we may construct the correspondence $E^*\odot F$ from $\cB$ to $\cC$ and get back $F$ as $F=E\odot(E^*\odot F)$. (Here $E$ need not be full.)

\brem \label{Motivation}
It is the last property in the preceding summary that furnishes the essential idea for the representation theory in Theorem \ref{strirepthm}: Given a strict unital homomorphism $\vt\colon \sB^a(E) \rightarrow \sB^a(F)$, we simply set $F_\vt=E^*\odot F$. The proof of Theorem \ref{strirepthm} (and the proof of Theorem \ref{assthm} for \nbd{W^*}modules) amounts to no more than working out the last part of the summary for explicit identifications. To understand the statement and proof, it is sufficient to know the tensor product of correspondences and the strict topology. It is \it{not} necessary to understand that $E$ is a Morita equivalence from $\sK(E)$ to $\cB_E$. This observation explains only why the construction works and how to find it. The reader familiar with tensor products and the strict topology could have started reading this note immediately with Theorem \ref{strirepthm}.

On the other hand, the preceding identifications are made only up to isomorphism. 
 The identification in the theorem is explicit and a good portion of this note is dedicated to showing that what intuition suggests 
 is compatible with the explicit identifications chosen in the theorem. In fact, we do not have much choice in writing down the identifications, because all objects are distinguished by universal properties. Taking into account these universal properties, all identifications become essentially unique; see [\refcite{Ske03p1}] for details.
\erem

\bitemp[Theorem (Representation Theorem).]\label{strirepthm}
Let $E$ be a Hilbert \nbd{\cB}module, let $F$ be a Hilbert \nbd{\cC}module and let $\vt\colon\sB^a(E)\rightarrow\sB^a(F)$ be a strict unital homomorphism. (In other words, $F$ is a correspondence from $\sB^a(E)$ to $\cC$ with strict left action and, thus, also a correspondence from $\sK(E)$ to $\cC$.) Then $F_\vt:=E^*\odot F$ is a correspondence from $\cB$ to $\cC$ and the formula
\beqn{
u(x_1\odot(x_2^*\odot y))
~:=~
\vt(x_1x_2^*)y
}\eeqn
defines a unitary
\beqn{
u
\colon
E\odot F_\vt
~\longrightarrow~
F
}\eeqn
such that
\beqn{
\vt(a)
~=~
u(a\odot\id_{F_\vt})u^*.
}\eeqn
\eitemp

\proof
To see that $u$ is isometric, simply compute:
\beqn{
\AB{x_1\odot(x_2^*\odot y),x'_1\odot({x'_2}^*\odot y')}
~=~
\AB{y,\vt(x_2\AB{x_1,x'_1}{x'_2}^*)y'}
~=~
\AB{\vt(x_1x_2^*)y,\vt(x'_1{x'_2}^*)y'}.
}\eeqn
Since $\sK(E)$ admits a bounded approximate unit $\bfam{u_\lambda}$ consisting of finite-rank operators and since $\vt$ is strict, we conclude that for all $y\in F$
\beqn{
y
~=~
\lim_\lambda\vt(u_\lambda)y
~\in~u(E\odot F_\vt).
}\eeqn
Thus, $u$ is also surjective. Finally,
\beqn{
\vt(a)u(x_1\odot x_2^*\odot y)
~=~
\vt(ax_1x_2^*)y
~=~
u(ax_1\odot x_2^*\odot y)
~=~u(a\odot\id_{F_\vt})(x_1\odot x_2^*\odot y),
}\eeqn
proving that $\vt(a)=u(a\odot\id_{F_\vt})u^*$ ($a\in\sB^a(E)$).~\qed

\bex\label{B(H)ex}
Returning to the normal (and, therefore, strict) representation $\vt$ of $\sB(G)$ on $H$ with which we began this note, we find that $H_\vt=G^*\odot H$, a tensor product that is balanced over $\sB(G)$, with the action on $H$ implemented by $\vt$.  It thus has much more structure than is revealed by thinking of it simply as a Hilbert space. 
 Of course, it is possible to define $H_\vt$ directly as the \it{Kolmogorov decomposition} of the positive definite kernel $\bfam{(g_1^*,h_1),(g_2^*,h_2)}\mapsto\AB{h_1,\vt(g_1g_2^*)h_2}$ on $G^*\times H$. However, then one must define a unitary $g'\otimes(g^*,h)\mapsto \vt(g'g^*)h$ and verify the desired properties. Like some other proofs based on the kernel approach, it appears largely unmotivated. The procedure does not explain what suggested the particular kernel, nor does it help with details.  In particular, it does not explain why the kernel actually \it{is} positive definite, nor does it explain why the resulting \it{Kolmogorov decomposition} has the desired properties.  Only the fact that the canonical map $G^* \times H \mapsto H_{\vt}$ is \it{balanced} over $\sB(G)$ reveals the structure of $H_{\vt}$.

On the other hand, while the internal structure of $H_{\vt} = G^*\odot H$ is revealed through its formulation as a balanced tensor product, it is difficult to capture $H_\vt$ in terms of concrete spaces. (After all, $G^*\odot H$ is an abstract space that must be constructed.) For example, in Bhat's approach [\refcite{Bha96}], it is possible to identify $H_\vt$ with the subspace $\vt(g_0g_0^*)H$ of $H$ via $g^*\odot h\mapsto\vt(g_0g^*)h$ for every unit vector $g_0$ in $G$. But there is no canonical choice for $g_0$ and the subspace, of course, depends on $g_0$. Arveson [\refcite{Arv89}], on the other hand, captures a space $H'_\vt$, isomorphic to $H_\vt$, as an intertwiner space. However, this intertwiner space, formally, should be viewed not so much as a Hilbert space but rather as a correspondence from $\C'$ to $\C'$ ($\C'$ is the commutant of $\C$ in $\sB(\C)=\C$). This may sound like a total triviality, but it becomes quite significant when applied to compositions, where we find the ``opposite'' (i.e.\ covariant) correspondence between endomorphisms and tensor products, as opposed to the natural, contravariant correspondence as described in Theorem \ref{assthm} below. (Tsireslon [\refcite{Tsi00p1}] showed that the product systems constructed by Bhat and by Arveson in that way, indeed, need not be isomorphic.) It should be noted that Bhat's use of the representation theorem for $\sB(H)$ inspired the second author's approach to the representation theorem for $\sB^a(E)$ under the assumption that $E$ has a unit vector, i.e., a vector $\xi$ such that $\AB{\xi,\xi} = \U$ [\refcite{Ske02}].  On the other hand, Arveson's analysis of endomorphisms inspired the second author's approach to the representation theorem for endomorphisms of $\sB^a(E)$, when $E$ is a full von Neumann module, in [\refcite{Ske03c}].  These two proofs are quite different in spirit and are different from the proof of Theorem \ref{strirepthm}.  An analysis of the differences will appear in [\refcite{Ske03p1}]; see also the survey [\refcite{Ske05a}].
\eex

\brem
Also, we should note that Rieffel's \it{imprimitivity theorem} [\refcite{Rie74}, Theorem 6.29] is Theorem \ref{strirepthm} when $F$ is a Hilbert space and with $\sF(E)$ instead of $\sB^a(E)$ (plus, the irrelevant technical modifications that $E$ be full and that $\cB$ may be only a pre-\nbd{C^*}algebra). In fact, our proof of Theorem \ref{strirepthm} is, really, a variation of the argument he employed for $\sF(E)$ and, then, extension to $\sB^a(E)$ and allowing a Hilbert module $F$ as representation space. We quote Rieffel's original result in Corollary \ref{Imprimitivity} and Rieffel's original argument is caputered as a part of the proof of Corollary \ref{NonUnital}.
\erem

\brem
When, in the setting of Theorem \ref{strirepthm} $\vt$ is not assumed to be unital, then $u$ is an isomorphism from $E\odot F_\vt$ onto the submodule $\vt(\U)F$ of $F$. This settles the complete treatment of strict representations (and, of course, of normal representations after Theorem \ref{normrepthm}) of $\sB^a(E)$ on another Hilbert module (\nbd{W^*}module) by simply applying Theorem \ref{strirepthm} to the unital representation $a\mapsto\vt(a)\upharpoonright\vt(\U)F$ of $\sB^a(E)$ on $\vt(\U)F$ and, then, embedding $\sB^a(\vt(\U)F)$ as the corner $\vt(\U)\sB^a(F)\vt(\U)$ in $\sB^a(F)$.
\erem

The following theorem shows that the choice of the multiplicity correspondence in the representation theorem is unique (provided fullness is assumed).

\bthm\label{unithm}
Let $E$ be a full Hilbert \nbd{\cB}module and let $F_1$ and $F_2$ be correspondences from $\cB$ to $\cC$. Suppose there exists a unitary $u\colon E\odot F_1\rightarrow E\odot F_2$ such that
\beqn{
u(a\odot\id_{F_1})
~=~
(a\odot\id_{F_2})u,
}\eeqn
i.e.\ $u$ is an \hl{isomorphism} of correspondences from $\sB^a(E)$ to $\cC$ (a \nbd{\sB^a(E)}\nbd{\cC}linear unitary). Then
\beq{\label{isoF}
F_1
~\longrightarrow~
\cB\odot F_1
~\longrightarrow~
E^*\odot E\odot F_1
~\xrightarrow{~~\id_{E^*}\odot u~~}~
E^*\odot E\odot F_2
~\longrightarrow~
\cB\odot F_2
~\longrightarrow~
F_2
}\eeq
is the unique isomorphism $u_{F_1,F_2}\colon F_1\rightarrow F_2$ such that $\id_E\odot u_{F_1,F_2}=u$.

In particular, if $u\colon E\odot F_1\rightarrow F=E\odot E^*\odot F$ is an isomorphism of strict correspondences from $\sB^a(E)$ to $\cC$ (i.e., in particular, if $F$ carries a strict unital representation $\vt$ of $\sB^a(E)$), then $u_{F_1,F_\vt}\colon F_1\rightarrow F_\vt$ is the unique isomorphism such that $\id_E\odot u_{F_1,F_\vt}=u$.
\ethm

\proof
The proof is immediate simply by recalling the canonical identifications of $F_i$ and $\cB\odot F_i$ and of $E^*\odot E$ and $\cB$ as correspondences.~\qed

\brem\label{nonfrem}
Actually, for the uniqueness part of the theorem the condition that $E$ be full is also necessary. Indeed, if the range ideal $\cB_E=E^*\odot E$ (viewed as correspondence from $\cB$ to $\cB$) is not $\cB$, then both $\cB$ and $\cB_E$ are correspondences from $\cB$ to $\cB$ inducing the identity representation, i.e.\ $E\odot\cB=E\odot\cB_E=E$ and $a\odot\id_\cB=a\odot\id_{\cB_E}=a$.
\erem

\bob\label{nonfob}
Also if $E$ is not necessarily full, the middle step in \eqref{isoF} is still an isomorphism $\cB_E\odot F_1=E^*\odot E\odot F_1\rightarrow E^*\odot E\odot F_2=\cB_E\odot F_2$. Clearly, the submodule $\cB_E\odot F_i=E^*\odot(E\odot F_i)$ of $F_i$ is the correspondence constructed as in Theorem \ref{strirepthm} from the representation $a\mapsto a\odot\id_{F_i}$ on $E\odot F_i$.

If we view $E$ as a full Hilbert \nbd{\cB_E}module, and if we view $\cB_E$ as a correspondence from $\cB_E$ to $\cB_E$ so that $\cB_E\odot F_i$ is a correspondence from $\cB_E$ to $\cC$, then under the assumptions of Theorem \ref{strirepthm} we obtain: $F_\vt$ is the unique correspondence from $\cB_E$ to $\cC$ such that $F=E\odot F_\vt$ (where the tensor product is taken over $\cB_E$) and $\vt(a)=a\odot\id_{F_\vt}$. If $F_2$ is a correspondence from $\cB$ to $\cC$ such that $F=E\odot F_2$ and $\vt(a)=a\odot\id_{F_2}$, in other words, if the hypothesis of Theorem \ref{unithm} are satisfied, with $F_1=F_\vt$, then $F_2\supset\cB_E\odot F_2\cong F_\vt$ via the isomorphism
\beq{\label{F1fu}
F_\vt
~\longrightarrow~
\cB_E\odot F_\vt
~\longrightarrow~
E^*\odot E\odot F_\vt
~\xrightarrow{~~\id_{E^*}\odot u~~}~
E^*\odot E\odot F_2
~\subset~
\cB\odot F_2
~\longrightarrow~
F_2
}\eeq
and $u$ is recovered as
\beq{\label{F1F2}
E\odot F_\vt\ni x\odot y\mapsto x\odot y\in E\odot F_2
}\eeq
where $y$ is, first, an element of $F_\vt$ that is, then, interpreted via \eqref{F1fu} as an element of $F_2\supset F_{\vt}$. Roughly speaking, the elements $x\odot y$ with $y\in\cB_EF_2$ are total in $E\odot F_2$.

Note that Remark \ref{nonfrem} does not contradict the uniqueness statement. $\cB$ is neither a Hilbert \nbd{\cB_E}module (its inner product does not take values in $\cB_E$) nor a correspondence from $\cB_E$ to $\cB$ (because the left action of $\cB_E$ is degenerate). We can reformulate it in the following way: The correspondence $F_\vt=E^*\odot F$ of Theorem \ref{strirepthm} is the unique correspondence from $\cB$ to $\cC$ (fulfilling the stated property) that is also a correspondence from $\cB_E$ to $\cC$, that is, $F_\vt=\cB_E\odot F_\vt$ (tensor product over $\cB$!).
\eob

\bcor\label{Morcor}
A full Hilbert \nbd{\cB}module $E$ and a full Hilbert \nbd{\cC}module $F$ have \hl{strictly isomorphic} operator algebras (the isomorphism and its inverse are strict mappings), if and only if there is a Morita equivalence $M$ from $\cB$ to $\cC$ such that $F\cong E\odot M$.
\ecor

\proof
Suppose $M$ is a Morita equivalence $M$ from $\cB$ to $\cC$. Then $(E\odot M)\odot M^*=E$ in the canonical identifications $(E\odot M)\odot M^*=E\odot(M\odot M^*)$, $M\odot M^*=\cB$ and $E\odot\cB=E$. In these identifications the strict homomorphisms $\sB^a(E)\ni a\mapsto a\odot\id_M$ and $\sB^a(E\odot M)\ni a\mapsto a\odot\id_{M^*}$ are inverses of each other. In particular, $\sB^a(E)$ and $\sB^a(E\odot M)$ are strictly isomorphic. If now also $F\cong E\odot M$, then $\sB^a(F)\cong\sB^a(E\odot M)$ (including the strict topology that is induced by $\sK(F)\cong\sK(E\odot M)$) so that, finally,  $\sB^a(F)$ and $\sB^a(E)$ are also strictly isomorphic.

Now suppose that $\vt\colon\sB^a(E)\rightarrow\sB^a(F)$ is a strict homomorphism with strict inverse. By Theorem \ref{strirepthm} there exist a correspondence $M$ from $\cB$ to $\cC$ and a correspondence $N$ from $\cC$ to $\cB$ such that $\vt$ is unitarily equivalent to $\sB^a(E)\ni a\mapsto a\odot\id_M$ and $\vt^{-1}$ is unitarily equivalent to $\sB^a(F)\ni a\mapsto a\odot\id_N$. Therefore, $\id_{\sB^a(E)}=\vt^{-1}\circ\vt$ is unitarily equivalent to $a\mapsto a\odot\id_M\odot\id_N=a\odot\id_{M\odot N}$. By Theorem \ref{unithm} this identifies (up to isomorphism) the correspondence $M\odot N$ as $\cB$. Similarly, the correspondence $N\odot M$ is isomorphic to $\cC$. In other words, $N\cong\cC\odot N\cong M^*\odot M\odot N\cong M^*\odot\cB\cong M^*$ and $M$ is a Morita equivalence.~\qed

\brem\label{Morrem}
Corollary \ref{Morcor} extends to the setting of nonfull Hilbert modules if we view them as Hilbert modules over their range ideals.  That is, if $E$ and $F$ are not necessarily full \nbd{\cB} and \nbd{\cC}modules, respectively, then Corollary \ref{Morcor} applies if we view $E$ and $F$ as full \nbd{\cB_E} and \nbd{\cC_F}modules, respectively. That is, if $\sB^a(E)$ and $\sB^a(F)$ are strictly isomorphic, then all we can conclude is that $\cB_E$ and $\cC_F$ are Morita equivalent.  On the other hand, a Morita equivalence between $\cB_E$ and $\cC_F$ induces an isomorphism between $\sB^a(E)$ and $\sB^a(F)$. (Note that for this it is not necessary that $\cB$ and $\cC$ are be Morita equivalent.)
\erem

\brem
We note that an isomorphism $\vt\colon\sB^a(E)\rightarrow\sB^a(F)$ is bistrict, if and only if both $\vt$ and $\vt^{-1}$ take the compacts into (and, therefore, onto) the compacts. (The ``if'' direction is obvious. The ``only if'' direction follows from the simple calculation $(x_1(m_1m_2^*)x_2^*)\odot\id_M=(x_1\odot m_1)(x_1\odot m_1)^*$, where $m_1m_2^*\in\sK(M)=\cB$. Therefore, choosing a bounded approximate unit $\bfam{u_\lambda}$ for $\cB=\sK(M)$ in the finite rank operators $\sF(M)$, we see that $x_1x_2^*\odot\id_M$ is approximated by finite-rank operators.) There are circumstances when the hypothesis of strictness in Corollary \ref{Morcor} is automatically fulfilled.  Anoussis and Todorov [\refcite{AnTo05}, Corollary 2.5] recently proved that if $\cB$ is unital and separable, and if $E$ and $F$ are countably generated over $\cB$, then any isomorphism from $\sB^a(E)$ onto $\sB^a(F)$ automatically maps $\sK(E)$ onto $\sK(F)$. A moment's reflection reveals that the same conclusion is valid when $E$ and $F$ are Hilbert modules over different \nbd{C^*}algebras \it{which need not be unital}. The way to remove the assumptions is to note that $E$ and $F$ are Hilbert modules over the bigger \nbd{C^*}algebra $(\cB\oplus\cC)^+$ (still separable and with an artificial new unit) in an obvious fashion. This, of course, does not affect any of the algebras $\sB^a(E)$, $\sK(E)$, $\sB^a(F)$, $\sK(F)$. 
\erem

Corollary \ref{Morcor} suggests the problem of identifying the relation between the pair of correspondences associated with two homomorphisms $\vt_1,\vt_2$ and the correspondence of their composition $\vt_2\circ\vt_1$. For the full case (implying uniqueness of the correspondences) we find, exactly as in the proof of Corollary \ref{Morcor}, that the latter is isomorphic to the tensor product (in contravariant order) of the former. However, these are identifications up to isomorphism and now we want to know whether our concrete identifications in Theorem \ref{strirepthm} are associative, i.e.\ we want to know whether they are compatible with the canonical identifications among multi-fold tensor products (making brackets superfluous) and other canonical identifications like $E^*\odot E=\cB$ via $x^*\odot x'=\AB{x,x'}$ and $E\odot E^*=\sK(E)$ via $x\odot x'^*=xx'^*$. (See the detailed discussion about the crucial difference between the statements ``{equal up to isomorphism}'' and ``{equal up to \hl{canonical} isomorphism}'' in [\refcite{Ske03p1}].) As an additional advantage, the results hold without the condition that the modules are full.

The generality in what follows is, of course, not necessary, but the notational convenience is considerable. All identifications are meant in the sense of Theorem \ref{strirepthm} ($F=E\odot F_\vt$ and $\vt(a)=a\odot\id_{F_\vt}$) and the other canonical identifications discussed above.

\bthm\label{assthm}
Let $E_i$ be Hilbert \nbd{\cB_i}modules ($i=1,2,\ldots$) and let $\vt_{i+1,i}\colon\sB^a(E_i)\rightarrow\sB^a(E_{i+1})$ be unital strict homomorphisms. For $1\le k<\ell$ set $\vt_{\ell,k}=\vt_{\ell,\ell-1}\circ\vt_{\ell-1,\ell-2}\circ\ldots\circ\vt_{k+1,k}$ and put $\vt_{i,i}=\id_{\sB^a(E_i)}$. For $1\le k\le\ell$ denote by $E_{\vt_{\ell,k}}$ the correspondence from $\cB_k$ to $\cB_\ell$ generating the homomorphism $\vt_{\ell,k}$ as in Theorem \ref{strirepthm}, while $E_{k,\ell}$ denotes the right Hilbert \nbd{\cB_\ell}module $E_\ell$ considered as a correspondence from $\sB^a(E_k)$ to $\cB_\ell$ with left action via $\vt_{\ell,k}$. Therefore, $E_{\vt_{\ell,k}}=E_k^*\odot E_{k,\ell}$, and $E_{k,\ell}\odot E_\ell^*$ is canonically isomorphic to the Hilbert \nbd{\sK(E_\ell)}module $\sK(E_\ell)$, when considered as correspondence from $\sK(E_k)$ to $\sK(E_\ell)$.

\begin{enumerate}
\item\label{1}
For every choice of $1\le i_1<\ldots<i_j$ we define an isomorphism
\beqn{
(E_{i_1,i_2}\odot E_{i_2}^*)\odot(E_{i_2,i_3}\odot E_{i_3}^*)\odot(E_{i_3,i_4}\odot E_{i_4}^*)\odot\ldots\odot(E_{i_{j-1},i_j}\odot E_{i_j}^*)
~\longrightarrow~
(E_{i_1,i_j}\odot E_{i_j}^*)
}\eeqn
by setting
\beqn{
a_2\odot a_3\odot a_4\odot\ldots\odot a_j
~\longmapsto~
(\ldots((a_2a_3)a_4)\ldots)a_j
}\eeqn
($a_2$ acts to the left on $a_3$, the result acts to the left on $a_4$, and so forth). Therefore, iteration of such isomorphisms is associative.

\item\label{2}
By tensoring with $E_{i_1}^*$ from the left and with $E_{i_j}$ ($=E_{i_j,i_j}$) from the right, we obtain an isomorphism
\beqn{
E_{\vt_{i_2,i_1}}\odot E_{\vt_{i_3,i_2}}\odot\ldots\odot E_{\vt_{i_j,i_{j-1}}}
~\longrightarrow~
E_{\vt_{i_j,i_1}}\odot E_{\vt_{i_j,i_j}}
~=~
E_{\vt_{i_j,i_1}}.
}\eeqn
Also these isomorphisms iterate associatively.
\end{enumerate}
\ethm

\proof
This follows simply by careful inspection of the identifications.~\qed

\brem
The isomorphism defined in \ref{1} extends from $E_{k,\ell}\odot E_\ell^*=\sK(E_\ell)$ (with left action of $\sK(E_k)$ via $\vt_{\ell,k}$) to the correspondences $\sB^a(E_\ell)$ from $\sB^a(E_k)$ to $\sB^a(E_\ell)$. (For instance, if $E_i=\cB$ ($\cB$ unital) for all $i$ so that $\vt_{i+1,i}$ are just unital endomorphisms of $\cB$, then $E_{\vt_{\ell,k}}$ and $E_{k,\ell}$ coincide and the isomorphisms from \ref{1} and \ref{2} are the same.)

Assertion \ref{2} suggests the notion of a nonstationary version of a product system. Indeed, if we set each $\cB_i=\cB$ and each $E_i=E$ and choose the $\vt_i$ to form a semigroup $\vt$ (indexed by $\N_0$ or $\R_+$), then by \ref{2} we obtain the product system $\vt$ as in [\refcite{Ske02,Ske03c}]. Passing from \ref{1} to \ref{2} is an operation of Morita equivalence for correspondences as considered in [\refcite{MuSo00}]. While the product system of \nbd{\sB^a(E)}modules in \ref{1} is \it{one-dimensional} and appears somewhat trivial, the Morita equivalent product system of \nbd{\cB}modules in \ref{2} can be very complicated. This appears already in the case when $E$ is a Hilbert space (which, of course, may be understood in an identical fashion), as studied by Arveson [\refcite{Arv89}].
\erem

We close this section with two results concerning \nbd{W^*}modules.  The first is an analogue of Theorem \ref{strirepthm} for \nbd{W^*}modules and normal homomorphisms.  The second generalizes  the well-known fact that a representation of $\sB(G)$ on a Hilbert space $H$ decomposes into the direct sum of two representations of $\sB(G)$.  One is normal (and, therefore, an amplification of the identity representation) and the other annihilates $\sK(G)$.  The representations of $\sB(G)$ that annihilate $\sK(G)$ are called \it{singular}.  Logically the second result should appear before the first one. But, there are a couple of technical explanations that would obstruct the simplicity of the argument as far as Theorem \ref{normrepthm} is concerned. 

For our purposes, the most convenient way to define a \hl{\nbd{W^*}module}, $E$ say, is to say that it is a \hl{self dual} Hilbert module over a \nbd{W^*}algebra $\cB$. ``Self duality'', in turn, means that every bounded right linear mapping $\Phi\colon E\rightarrow\cB$ has the form $\Phi x=\AB{y,x}$ for some (unique) $y\in E$.  If $E$ is a \nbd{W^*}module over a \nbd{W^*}algebra $\cB$, then $\sB^a(E)$ is also a \nbd{W^*}algebra. A correspondence over \nbd{W^*}algebras is a \hl{\nbd{W^*}correspondence}, if it is a right \nbd{W^*}module such that the canonical homomorphism determined by the left action is a normal mapping into $\sB^a(E)$. The \nbd{W^*}module $E$ is \hl{strongly full}, if the range of its inner product generates $\cB$ as a \nbd{W^*}algebra. Usually, we will leave out the word ``strongly''.

Every (pre-)Hilbert module over a \nbd{W^*}algebra admits a unique minimal self-dual extension (see Paschke [\refcite{Pas73}] or Rieffel [\refcite{Rie74a}]). By the tensor product of two \nbd{W^*}correspondences $E$ and $F$ we understand the self-dual extension of $E\odot F$ which we denote by $E\sbars{\odot}F$. Also, every (bounded) adjointable operator (like $a\odot\id_F$) on a (pre-)Hilbert module (like $E\odot F$) extends to a unique adjointable operator on the self-dual extension. The tensor product of \nbd{W^*}correspondences is a \nbd{W^*}correspondence.

\bthm\label{normrepthm}
Let $E$ be a \nbd{W^*}module over $\cB$, let $F$ be a \nbd{W^*}module over $\cC$ and let $\vt\colon\sB^a(E)$ $\rightarrow\sB^a(F)$ be a unital normal homomorphism. In other words, $F$ is a \nbd{W^*}correspondence from $\sB^a(E)$ to $\cC$. Then $F_\vt:=E^*\sbars{\odot}F$ is a \nbd{W^*}correspondence from $\cB$ to $\cC$ and the mapping
\beqn{
u(x_1\odot(x_2^*\odot y))
~:=~
\vt(x_1x_2^*)y
}\eeqn
defines a unitary
\beqn{
u
\colon
E\sbars{\odot}F_\vt
~\longrightarrow~
F
}\eeqn
such that
\beqn{
\vt(a)
~=~
u(a\odot\id_{F_\vt})u^*.
}\eeqn
\ethm

We leave it to the reader to state and prove the \nbd{W^*}versions of all the other results we proved for Hilbert modules (Theorem \ref{unithm}, Observation \ref{nonfob}, Corollary \ref{Morcor} and Theorem \ref{assthm}).

The following theorem is something of a mixture of the \nbd{C^*} and \nbd{W^*}frameworks. It deals with (arbitrary) representations of $\sB^a(E)$ in $\sB^a(F)$ where $E$ is a Hilbert module but $F$ is a \nbd{W^*}module. The critical properties of \nbd{W^*}modules, that enter the proof, are as follows. If $F_0$ is a pre-Hilbert submodule of $F$, then the \hl{bicomplement} $F_0^{\perp\perp}$ of $F_0$ in $F$ is a \nbd{W^*}module. In fact, it is the \nbd{W^*}module generated by $F_0$. Of course, $F=F_0^{\perp\perp}\oplus F_0^\perp$. (Of course, all this holds if $F_0$ is an arbitrary subset of $F$.) Moreover, an arbitrary bounded right-linear operator $a\colon F_0\rightarrow F$ extends uniquely to an element of $\sB^a(F)$ that vanishes on $F_0^\perp$.

\brem
These statements (as well as others needed for Theorem \ref{normrepthm}) are slightly tricky to prove in the setting of abstract \nbd{W^*}mod\-ules. They are more easily understood in the equivalent category of von Neumann modules; see Skeide [\refcite{Ske00b,Ske01}] and also Remark \ref{Rierem} below.
\erem

These are the major ingredients of a proof of the following theorem. We leave it to the reader to fill in the details.

\bthm\label{normalsingularrepthm}
Let $E$ be a Hilbert module over the \nbd{C^*}algebra $\cB$, let $F$ be a \nbd{W^*}module over the \nbd{W^*}algebra $\cC$ and let $\vt\colon \sB^a(E) \rightarrow \sB^a(F)$ be a representation.  Set $F_{sing}:=(\vt (\sK(E))F)^\perp$ and set $F_{ampl}:=F_{sing}^\perp=(\vt (\sK(E))F)^{\perp \perp}$.  Then:

\begin{enumerate}
\item

$F = F_{ampl} \oplus F_{sing}$.

\item
$\vt(a)F_{sing}=0$ for all $a\in \sK(E)$ so that, in particular, $F_{sing}$ is invariant under $\vt(\sK(E))$.

\item
$\vt(\sK(E)) F_{ampl}~=~\vt(\sK(E))F~\subseteq~F_{ampl}$.

\item
If we identify $E\odot E^* \odot F$ with the norm closure of $\ls\vt(\sK(E))F$ in $F$ via the map $x\odot y^*\odot z \mapsto \vt(xy*)z$, then for all $a \in \sB^a(E)$, the restriction of $\vt(a)$ to $F_{ampl}$ is the unique extension of $a\odot \id_{E^*\odot F}$ to $F_{ampl}$.
\end{enumerate}
\ethm

Note that a Hilbert space is a \nbd{W^*}module (over $\C$) so that Theorem \ref{normalsingularrepthm} applies, in particular, to representations of $\sB^a(E)$ on Hilbert spaces.

\bprop
The restriction of $\vt$ to $F_{ampl}$, $\vt\upharpoonright F_{ampl}$, is strict if and only if the norm closure of $\vt(EE^*)F$ in $F$ is $F_{ampl}$, that is, if and only if $E\odot E^*\odot F=\ol{E\odot E^*\odot F}^s$ where $\ol{\rule{0mm}{2mm}~~}^s$ indicates the unique self-dual extension.
\eprop

\proof
If that representation is strict, then our Theorem \ref{strirepthm} tells us that we recover $F_{ampl}$ in that way. On the other hand, if \nbd{C^*}tensor products are sufficient to recover $F_{ampl}$, then the representation is strict.~\qed

\lf
When $E$ is a Hilbert module over a (pre-)\nbd{\cB}algebra $\cB$, then $\sF(E)$ is Rieffel's imprimitivity algebra [\refcite{Rie74}, Definition 6.4]. His \it{imprimitivity theorem} [\refcite{Rie74}, Theorem 6.29] (that is, Corollary \ref{Imprimitivity}), actually, is about inducing representations of $\sF(E)$ on Hilbert spaces from (bounded) representations of $\cB$ on Hilbert spaces via $E$.  The following corollary (or better, its proof) captures most clearly Rieffel's ideas that lead to the representation theorem for $\sF(E)$ and generalizes them to representations on Hilbert or \nbd{W^*}modules. But, it captures also the core of what we contributed in order to extend that result to $\sB^a(E)$ in a condensed way.

\bcor\label{NonUnital}
Let $E$ be a pre-Hilbert module over a \nbd{C^*}algebra (as everywhere in these notes) or Hilbert module over a pre-\nbd{C^*}algebra (as in Rieffel [\refcite{Rie74}]). Every representation $\vt$ of the finite-rank operators $\sF(E)$ on a Hilbert module (a \nbd{W^*}module) $F$ that is generated by $\vt(EE^*)F$ as a Hilbert module (as a \nbd{W^*}module) extends to a unique representation of $\sB^a(E)$. (In particular, $\vt$ is bounded automatically and strict if it generates $F$ as a Hilbert module.) In the case of Theorem \ref{normalsingularrepthm}, $\vt\upharpoonright F_{ampl}$ is that unique extension of the representation $xy^*\mapsto\vt(xy^*)\upharpoonright\ls\vt(EE^*)F$ of $\sF(E)$.
\ecor

\proof
Every Hilbert module over a pre-\nbd{C^*}algebra may be considered also as a Hilbert module over the completion of that algebra. So we cover both cases, if we suppose that $E$ is a pre-Hilbert module over a \nbd{C^*}algebra.

Let $F$ be generated by $\vt(EE^*)F$ as a Hilbert module. Then the proof of Theorem \ref{strirepthm} (as in [\refcite{Rie74}]) shows that $F=E\odot E^*\odot F$ and $\vt(xy^*)=(xy^*)\odot\id_{E^*\odot F}$. The fact that $E$ might not be complete is irrelevant as the tensor product is assumed completed. The only critical question is whether the completion of the algebraic tensor product $E^*\ul{\odot}F$ is still a correspondence from $\cB$ to $\cC$ and this follows, because all representations of a \nbd{C^*}algebra by possibly unbounded operators on a pre-Hilbert module, acutally, are by bounded operators. (Note that this also shows that $E\odot E^*\odot F=\ol{E}\odot\ol{E}^{\,*}\odot F$.)

Clearly, if $\vt$ possesses an extension to $\sB^a(E)$, then this extension is uniquely determined as $a\odot\id_{E^*\odot F}$, because $\vt(a)\vt(xy^*)z=\vt(axy^*)z=(ax)\odot(y\odot z)$. In particular, this representation extends further to $\sB^a(\ol{E})$ and it is strict.

The extension result when $F$ a \nbd{W^*}module and generated as such by its submodule $E\odot E^*\odot F$ follows by the mentioned extension result for operators defined on generating submodules.~\qed

\brem
If $E$ is a pre-Hilbert module over a \nbd{W^*}algebra $\cB$ and if $\ol{E}^s$ is its unique extension as a \nbd{W^*}module, then it may be possible to extend $\vt$ to some portion of $\sB^a(\ol{E}^s)\supset\sB^a(\ol{E})$, but perhaps not to all of $\sB^a(\ol{E}^s)$.  We do not know what the precise details are.  However, $\vt$ extends to all of $\sB^a(\ol{E}^s)$ if and only if the left action of $\cB$ on $E^*\sodots F$ makes it a \nbd{W^*}correspondence.
\erem

\bitemp[{Corollary \protect[\refcite{Rie74}, Theorem 6.29].}]\label{Imprimitivity}
If $E$ is a full Hilbert \nbd{\cB}module over a (pre-) \nbd{C^*}al\-ge\-bra $\cB$ and if $\vt$ is a nondegenerate representation of $\sF(E)$ on a Hilbert space $H$, then $\vt$ is unitarily equivalent to a representation of $\sF(E)$ that is \hl{induced} via $E$ in the sense of Rieffel by a representation of $\cB$ (that is, a representation of the form $a\mapsto a\odot\id_G$ for some nondegenerate representation of $\cB$ on a Hilbert space $G$).
\eitemp

\section{Functors}

In this section we discuss the functorial aspects of tensoring modules with a fixed correspondence. We utilize our representation theorem to furnish a new proof of Blecher's \it{Eilenberg-Watts theorem} [\refcite{Ble97}]. We also discuss relations between his approach and ours.

We denote by $\eC^*_\cB$ the category whose objects are Hilbert \nbd{\cB}modules and whose morphisms are adjointable mappings.   Also, we denote by $_\cA\smash{\eC^*}\!\!\!_\cB$ the category whose objects are correspondences from $\cA$ to $\cB$ and whose morphisms are adjointable mappings that intertwine the actions of $\cA$ as morphisms. We denote the corresponding categories of \nbd{W^*}modules and \nbd{W^*}correspondences by $\eW^*_\cB$ and $_\cA\smash{\eW^*}\!\!\!_\cB$, respectively.  It is easily verified that $\eC^*_\cB$ and $_\cA\smash{\eC^*}\!\!\!_\cB$ are \nbd{C^*}categories in the sense of Ghez, Lima and Roberts [\refcite{GLR85}], while $\eW^*_\cB$ and $_\cA\smash{\eW^*}\!\!\!_\cB$ are \nbd{W^*}categories in their sense. We concentrate on the \nbd{C^*}categories and leave the simpler \nbd{W^*}categories to the reader.

There are two possibilities for tensoring with a fixed correspondence, namely, tensoring from the left and tensoring from the right. The second possibility is the one relevant for Blecher's \it{Eilenberg-Watts} theorem. But before we focus on this second possibility, we quickly review some well-known facts about the first possibility, which is related to Rieffel's \it{Eilenberg-Watts theorem} [\refcite{Rie74a}] and to which we return in Remark \ref{Rierem}. Let $E$ be a correspondence from $\cA$ to $\cB$. Then, for every \nbd{C^*}algebra $\cC$ we obtain a functor $\el_E^\cC\colon {_\cB\smash{\eC^*}\!\!\!_\cC}\rightarrow {~_\cA\smash{\eC^*}\!\!\!_\cC}$ by setting
\beqn{
\el_E^\cC(F)
~=~
E\odot F,
~~~
\el_E^\cC(a)
~=~
\id_E\odot a,
~~~~~~(F,F_1,F_2\in {~_\cB\smash{\eC^*}\!\!\!_\cC}; a\in\sB^{a,bil}(F_1,F_2)).
}\eeqn
It is also easy to see that  $\el_{E_1}^\cC\circ\el_{E_2}^\cC=\el_{E_1\odot E_2}^\cC$. If $M$ is a Morita equivalence from $\cA$ to $\cB$, then the functor $\el_M^\cC$ is an equivalence from $_\cB\smash{\eC^*}\!\!\!_\cC$ to $_\cA\smash{\eC^*}\!\!\!_\cC$. In particular, putting $\el_E:=\el_E^\C$, the functor $\el_M$ is an equivalence from the category of nondegenerate representations of $\cB$ (on Hilbert spaces) to that of representations of $\cA$ (on Hilbert spaces). Rieffel [\refcite{Rie74a}] calls $\cB$ and $\cA$ \hl{strongly Morita equivalent}, if such an $\el_M$ exists. Rieffel's \it{Eilenberg-Watts theorem} asserts that every (sufficiently regular) functor between the categories of representations of two \nbd{W^*}algebras $\cB$ and $\cA$ has the form $\el_E$ for a suitable \nbd{W^*}correspondence $E$ from $\cA$ to $\cB$.

Now let us focus on the second possibility --- tensoring on the right. Fix a correspondence $F$ from $\cB$ to $\cC$. Then we define a functor $\er_F\colon\eC^*_\cB\rightarrow\eC^*_\cC$ by setting
\beqn{
\er_F(E)
~=~
E\odot F,
~~~
\er_F(a)
~=~
a\odot\id_F,
~~~~~~(E,E_1,E_2\in\eC^*_\cB; a\in\sB^a(E_1,E_2)).
}\eeqn
$\er_F$ is a \hl{\nbd{*}functor} in the sense that $\er_F(a^*)=\er_F(a)^*$. It is also \hl{strict} in the sense that it is strictly continuous on bounded subsets of $\sB^a(E_1,E_2)\subset\sB^a(E_1\oplus E_2)$. (Note that for any \nbd{*}functor $\er$ this is equivalent to saying that $\er$ is strongly continuous on bounded subsets, because being strict means being \nbd{*}strongly continuous on bounded subsets, and being strongly continuous on bounded subsets, for a \nbd{*}functor, implies that the functor is \nbd{*}strongly continuous on bounded subsets.) Also the functor $\er_F$ is an equivalence (of categories), if and only if $F$ is a Morita equivalence. Further, we note that the functors $\er_F$ compose contravariantly to tensoring, i.e.\ $\er_{F_1}\circ\er_{F_2}=\er_{F_2\odot F_1}$.

Our goal is to show that every strict \nbd{*}functor $\er\colon\eC^*_\cB\rightarrow\eC^*_\cC$ arises as an $\er_F$ for a suitable correspondence $F$ from $\cB$ to $\cC$. More precisely, we show that $\er$ is \hl{naturally equivalent} to $\er_F$, for some $F$, meaning that for every object $E\in\eC^*_\cB$ there is an isomorphism $v_E\colon E\odot F\rightarrow\er(E)$ such that $\er(a)=v_{E_2}(a\odot\id_F)v_{E_1}^*$ for all objects $E_1,E_2\in\eC^*_\cB$ and all morphisms $a\in\sB^a(E_1,E_2)$.

If we focus on a single object $E\in\eC^*_\cB$, the restriction of $\er$ to $\sB^a(E)$ is a unital strict representation of $\sB^a(E)$ on $\er(E)$. Theorem \ref{strirepthm} then provides us with the correspondence $F_E:=E^*\odot\er(E)$ and an isomorphism $u_E\colon E\odot F_E \rightarrow \er(E)$ such that $\er(a)=u_E(a\odot\id_{F_E})u_E^*$ for all $a\in\sB^a(E)$. So, in order to prove the \it{Eilenberg-Watts theorem} we have to face two problems. First, we have to eliminate the dependence of $F_E$ on $E$.  That is, we want to find a single $F$  that works for all $E$. (Observation \ref{nonfob} will provide us with the necessary tools.) Second, we have to show that not only does the equation $\er(a)=v_E(a\odot\id_F)v_E^*$ hold for all $a\in\sB^a(E)$ (as implied by Theorem \ref{strirepthm}) but also the equation $\er(a)=v_{E_2}(a\odot\id_F)v_{E_1}^*$ must be satisfied for $a\in\sB^a(E_1,E_2)$ for every \it{pair} of Hilbert modules $E_1$ and $E_2$. In the sequel, we will achieve both goals by showing that, in a suitable sense, $\er$ behaves ``nicely'' with respect to direct sums.  That is, we show that a \nbd{*}functor $\er$ from $\eC^*_\cB$ to $\eC^*_\cC$ must be additive.

\bob
Since $\er$ is a functor, it can only be applied to objects $E$ and to morphisms $a$. It cannot be applied to elements $x\in E$. However, it can be applied to rank-one operators $xy^*$ $(x\in E_1,y\in E_2)$. If we do so, then we have factorizations like $\er(axy^*)=\er(a)\er(xy^*)$ or $\er(x(ay)^*)=\er(xy^*)\er(a^*)$. This key observation is central to our analysis and is used over and over again as a substitute for the (non-existent) values of $\er$ at points of $E$.
\eob

We start by listing several self-evident properties that every \nbd{*}functor from $\eC^*_\cB$ to $\eC^*_\cC$ must satisfy.

\bprop\label{*funprop}
Every \nbd{*}functor $\er\colon\eC^*_\cB\rightarrow\eC^*_\cC$ sends projections $p\in\sB^a(E)$ to projections $\er(p)\in\sB^a(\er(E))$ and, therefore, it sends partial isometries $w\in\sB^a(E_1,E_2)$ to partial isometries $\er(w)\in\sB^a(\er(E_1),\er(E_2))$.

Further, since $\er(\id_E)=\id_{\er(E)}$, the functor $\er$ sends isometries $v\in\sB^a(E_1,E_2)$ to isometries $\er(v)\in\sB^a(\er(E_1),\er(E_2))$ and, therefore, $\er$ sends unitaries $u\in\sB^a(E_1,E_2)$ to unitaries $\er(u)\in\sB^a(\er(E_1),\er(E_2))$.
\eprop

\brem
Recall that isometries between Hilbert modules need not be adjointable. In fact, an isometry is adjointable if and only if there exists a projection onto its range, i.e.\ if and only if its range is complemented; see, for instance, [\refcite{Ske01}, Proposition 1.5.13]. Fortunately, in direct sums there exist projections onto the direct summands and, therefore, the canonical injections are adjointable.
\erem

Let $E_i$ $(i=1,2)$ be objects in $\eC^*_\cB$ and denote their direct sum by $E:=E_1\oplus E_2$. Denote by $k_i\in\sB^a(E_i,E)$ the canonical injections $E_i\rightarrow E$ and by $p_i=k_ik_i^*\in\sB^a(E)$ the projections onto $E_i$. By Proposition \ref{*funprop} $\er(k_i)$ is an isometry and $\er(p_i)$ is a projection onto its range. Consequently $\er(k_i)\er(E_i)=\er(p_i)\er(E)$. We obtain the following chain of unitaries as a result.
\bmu{\label{EE_iiso}
E_i\odot F_E
~\xrightarrow{~~k_i\odot\id_{F_E}~~}~
(k_i\odot\id_{F_E})(E_i\odot F_E)
~=~
(p_i\odot\id_{F_E})(E\odot F_E)
\\
~\xrightarrow{~~u_E~~}~
\er(p_i)\er(E)
~=~
\er(k_i)\er(E_i)
~\xrightarrow{~~\er(k_i)^*~~}~
\er(E_i)
~\xrightarrow{~~u_{E_i}^*~~}~
E_i\odot F_{E_i}.
}\emu
The trick is to translate the part that contains the operator $k_i\odot\id_{F_E}$, for which we do not yet know any explicit relation with the operator $\er(k_i)$, into an expression involving the operator $p_i\odot\id_{F_E}$, which we know is unitarily equivalent to $\er(p_i)$. Clearly, the resulting unitary $E_i\odot F_E\rightarrow E_i\odot F_{E_i}$  intertwines the canonical actions $a\odot\id_{F_E}$ and $a\odot\id_{F_{E_i}}$ of $a\in\sB^a(E_i)$. (To see this we just have to observe that the relevant action of $a$ on an element $k_ix_i$ $(x_i\in E_i)$ is $(k_iak_i^*)(k_ix_i)=k_i(ax_i)$ and on an element $\er(k_i)z_i$ $(z_i\in\er(E_i))$ is $\er(k_iak_i^*)\er(k_i)z_i=\er(k_i)\er(a)z_i$.)

By Observation \ref{nonfob}, $F_{E_i}$ is isomorphic to the submodule $\cls\cB_{E_i}F_E$ of $F_E$. To fix the isomorphism explicitly, we read \eqref{EE_iiso} backwards. So let us choose an element $x_i\odot(y_i^*\odot z_i)$ in $E_i\odot F_{E_i}=E_i\odot(E_i^*\odot\er(E_i))$ and apply $(k_i^*\odot\id_{F_E})u_E^*\er(k_i)u_{E_i}$. We find
\bmun{
x_i\odot y_i^*\odot z_i
~\overset{u_{E_i}}{\longmapsto}~
\er(x_iy_i^*)z_i
~\overset{\er(k_i)}{\longmapsto}~
\er(k_i)\er(x_iy_i^*)z_i
\\[1ex]
~=~
\er((k_ix_i)y_i^*)z_i
~=~
\er((k_ix_i)y_i^*)\er(k_i^*)\er(k_i)z_i
~=~
\er((k_ix_i)(k_iy_i)^*)(\er(k_i)z_i)
\\
~\overset{u_E^*}{\longmapsto}~
k_ix_i\odot(k_iy_i)^*\odot\er(k_i)z_i
~\overset{k_i^*\odot\id_{F_E}}{\longmapsto}~
x_i\odot(k_iy_i)^*\odot\er(k_i)z_i.
}\emun
This establishes that $y_i^*\odot z_i\mapsto(k_iy_i)^*\odot\er(k_i)z_i$ is the unique canonical embedding $F_{E_i}\rightarrow F_E$ discussed in Observation \ref{nonfob} (generalizing \eqref{isoF} to the not--necessarily--full case).

We apply this discussion to the direct sum $\cB\oplus E$,  where $E$ is any object in $\eC^*_\cB$, and where $k_\cB$ and $k_E$ are the canonical injections. Consider the direct summand $\cB$. Since $\cB$ is full, we have $F_\cB\cong F_{\cB\oplus E}$ and $b^*\odot z\mapsto(k_\cB b)^*\odot\er(k_\cB)z$ $(b\in\cB,z\in\er(\cB))$ is the isomorphism. On the other hand, if we focus on the direct summand $E$ we see that $F_E\cong\cls\cB_EF_{\cB\oplus E}$ and $y^*\odot z_E\mapsto(k_Ey)^*\odot\er(k_E)z_E$ $(y\in E,z_E\in\er(E))$ is the isomorphism. Moreover, we know that $E\odot F_{\cB\oplus E}\cong E\odot F_E$. The former is a subset $k_EE\odot F_{\cB\oplus E}$ of $(\cB\oplus E)\odot F_{\cB\oplus E}$ and the latter is isomorphic to $\er(E)$. Setting $F:=F_\cB$, we find that
\bmun{
x\odot(b^*\odot z)
~\longmapsto~
k_Ex\odot(k_\cB b)^*\odot\er(k_\cB)z
~\longmapsto~
\er((k_Ex)(k_\cB b)^*)\er(k_\cB)z
\\
~=~
\er(k_E)\er(x(k_\cB b)^*)\er(k_\cB)z
~\longmapsto~
\er(x(k_\cB b)^*)\er(k_\cB)z
~=~
\er(xb^*)\er(k_\cB)^*\er(k_\cB)z
~=~
\er(xb^*)z
}\emun
defines a unitary $v_E\colon E\odot F\rightarrow\er(E)$.

\brem
Of course, we could define $v_E$ immediately via the formula $x\odot(b^*\odot z)\mapsto\er(xb^*)z$. Is is easy to check that this defines an isometry. However, only the above chain of mappings together with the explanations preceding them show that the isometry is onto $\er(E)$.
\erem

\brem\label{F=r(B)rem}
As a Hilbert \nbd{\cC}module, $F$ coincides with $\er(\cB)$ and the left action of $b\in\cB$ is $\er(b)$ when $b$ is considered as an element of $\sB^a(\cB)\supset\cB$. The element $b^*\odot z$ of $F$ corresponds to the element $\er(b)^*z\in\er(\cB)$.
\erem

To conclude the argument we must show that the map $E\mapsto v_E$ is a natural transformation. Choose $a\in\sB^a(E_1,E_2)$ and $x\in E_1,b\in\cB,z\in\er(\cB)$. Then
\bmun{
v_{E_2}(a\odot\id_F)(x\odot(b^*\odot z))
~=~
v_{E_2}(ax\odot(b^*\odot z))
\\
~=~
\er(axb^*)z
~=~
\er(a)\er(xb^*)z
~=~
\er(a)v_{E_1}(x\odot(b^*\odot z)).
}\emun

We have, thus, proved the following.

\bitemp[Theorem (Eilenberg-Watts theorem).]\label{EWthm}
Let $\er\colon\eC^*_\cB\rightarrow\eC^*_\cC$ be a strict
\nbd{*}functor. Then $F=\cB^*\odot\er(\cB)$ is a correspondence in $_\cB\smash{\eC^*}\!\!\!_\cC$ such that $\er_F$ is naturally equivalent to $\er$ via the natural transformation $v_E\colon\er_F(E)\rightarrow\er(E)$ defined by setting $v_E(x\odot(b^*\odot z))=\er(xb^*)z$.

Moreover, $F$ is unique in $_\cB\smash{\eC^*}\!\!\!_\cC$. That is, if $\hat{F}\in_\cB\smash{\eC^*}\!\!\!_\cC$ is another correspondence such that $\er_{\hat{F}}$ is naturally equivalent to $\er$, then $F\cong\hat{F}$.
\eitemp

\bcor
The functor $\er$ is determined by its restriction to the object $\cB$ or, more generally, to any full Hilbert \nbd{\cB}module $E$, and its endomorphisms.
\ecor

We complement this with some more consequences of Theorems \ref{strirepthm} and \ref{unithm}.

\bprop
Every strict unital representation $\vt$ of $\sB^a(E)$ on a Hilbert \nbd{\cC}module $F$ extends to a strict \nbd{*}functor $\er_{F_\vt}$. This functor is unique (up to natural equivalence), if $E$ is full.
\eprop

We list a few properties that hold for $\er_F$ and, therefore, for every strict \nbd{*}functor.  The first assertion actually is true for every (not-necessarily strict) \nbd{*}functor [\refcite{GLR85}, Page 84].

\bprop\label{r_Fprop}
Let $F$ be a correspondence from $\cB$ to $\cC$. Then:
\begin{enumerate}
\item
$\er_F$ is \hl{contractive}, i.e.\ $\norm{\er(a)}\le\norm{a}$ for every morphism.

\item\label{Pii}
$\er_F$ extends to a \hl{strong} contractive functor $\eC_\cB\rightarrow\eC_\cC$ between the categories with the same objects but with bounded module maps as morphisms. (Here strong means that every restriction $\er_F\upharpoonright\sB^r(E_1,E_2)$ is strongly continuous on bounded subsets.)

\item
This extension \hl{respects embeddings}, i.e.\ if $v\in\sB^r(E_1,E_2)$ is an isometry, then $\er_F(v)$ is an isometry, too, and $\er_F(v)\er_F(E_1)=\er_F(vE_1)\subset\er_F(E_2)$.
\end{enumerate}
\eprop

\brem
Theorem \ref{EWthm} holds for an arbitrary \hl{full} subcategory of $\eC_\cB^*$, provided that the subcategory contains at least one full Hilbert \nbd{\cB}module $E_f$ and with every object $E$ also the direct sum $E_f\oplus E$. (Recall that to say a subcategory $\euf{D}$ of a category $\eC$ is full means that given any two objects in $\euf{D}$, all the morphisms between them that appear in $\eC$ appear also in $\euf{D}$.)
\erem

\brem%
Theorem \ref{strirepthm} plays a central role in the proof of Theorem \ref{EWthm}, and one may wonder if it is possible somehow to obtain Theorem \ref{strirepthm} from Theorem \ref{EWthm}.  The problem is that given an $E$ in $\eC^*_{\cB}$ and a homomorphism $\vt\colon \sB^a(E) \rightarrow \sB^a(F)$ for some $F$ in $\eC^*_{\cC}$, then to apply Theorem \ref{EWthm} to $\vt$, one must figure out how to extend $\vt$ directly to a functor $\er_{\vt}$ from $\eC^*_{\cB}$ to $\eC^*_{\cC}$ without passing through our Theorem \ref{strirepthm}.  (Blecher's argument shows that it suffices, really, to extend $\vt$ to the full subcategory of $\eC^*_{\cB}$ generated by the objects $\cB,E,\cB\oplus E$.) For then we could follow Blecher and apply $\er_{\vt}$ to $\cB$ to obtain a correspondence $F_{\cB}$ from $\cB$ to $\cC$ such that $F=E\odot F_{\cB}$.  Whether or not there is any evident functorial extension $\er_{\vt}$ of $\vt$ we don't know, but in any case, the process of passing through Theorem \ref{EWthm} to prove Theorem \ref{strirepthm} seems far more complicated than the direct argument we provided for Theorem \ref{strirepthm}.
\label{nbackrem}
\erem

\brem
We wish to mention also that the hypotheses in Blecher's version of the \it{Eilen\-berg-Watts theorem} are slightly different from ours. He considers functors between categories of Hilbert modules where \it{all} the bounded right module mappings are morphisms (as in Proposition \ref{r_Fprop}\eqref{Pii}). He must assume that his functors are bounded.  When one of his functors is restricted to our category, the restriction is by definition a \nbd{*}functor.
\erem

\brem\label{Rierem}
Finally, we would like to compare briefly Rieffel's \it{Eilenberg-Watts theorem} [\refcite{Rie74a}] and Blecher's [\refcite{Ble97}, Theorem 5.4].  For a full account, see [\refcite{Ske05p}] (and [\refcite{Ske03p1}]).  First of all, Rieffel's result is formulated for categories of normal unital representations of \nbd{W^*}algebras.  Consequently, when thinking about Blecher's version, one must reformulate it in terms of \nbd{W^*}categories.  This is not difficult and, in fact, because all bounded linear operators on a \nbd{W^*}module are automatically adjointable, such a reformulation appears more familiar. Rieffel shows that a normal \nbd{*}functor $\el$ between categories of normal unital representations of \nbd{W^*}algebras is naturally equivalent to a functor $\el_F$ for a suitable \nbd{W^*}correspondence $F$. To connect this with Blecher's representation of a \nbd{*}functor in terms of tensoring on the right, it is convenient to reformulate the whole discussion of \nbd{W^*}categories in the context of von Neumann categories, where, recall, one deals with concrete von Neumann algebras acting on specific Hilbert spaces and with modules that are realized as concrete spaces of operators, possibly between different Hilbert spaces.  When this perspective is adopted, one finds a natural duality which carries a correspondence $F$ from a von Neumann algebra $\cB$ to a von Neumann algebra $\cC$ to a correspondence, denoted $\f (F)$ or $F'$, from the \it{commutant} of $\cC$, $\cC'$, to the \it{commutant} of $\cB$, $\cB'$. This process is indeed a normal \nbd{*}functor that generalizes naturally the process of forming the commutant of a von Neumann algebra, and this explains the notation and terminology: $\f (F) = F'$ is called the \hl{commutant} of $F$. (The case $\cC=\cB$ was discussed in [\refcite{Ske03c}] and, independently, in [\refcite{MuSo04}], and generalized to different algebras in [\refcite{MuSo05}].) The restriction to $\cC=\C$ gives a duality, still denoted $\f$ that is implicit in [\refcite{Rie74a}] between the category of concrete von Neumann \nbd{\cB}modules and the category of normal unital representations of $\cB'$ on Hilbert spaces. Under $\f$ a functor $\er$ is transformed into the functor $\el=\f\circ\er\circ\f$, with $\f\circ\el\circ\f=\er$ since $\f\circ\f$ is the identity functor. Moreover, a functor $\er_F$ for a von Neumann correspondence $F$ from $\cB$ to $\cC$ translates into the functor $\el_{F'}=\f\circ\er_F\circ\f$ where $F'$ is the commutant of $F$. Thus, the duality between von Neumann correspondences and their commutants provides also a duality between Rieffel's and Blecher's \it{Eilenberg-Watts theorems}. Consequently, the two results are logically equivalent and any proof for one gives, by duality, a proof for the other.
\erem

\setlength{\baselineskip}{2.5ex}


\begin{thebibliography}{10}

\bibitem{AnTo05}
M.~Anoussis and I.~Todorov, \emph{{Compact operators on Hilbert modules}},
  Proc.\ Amer.\ Math.\ Soc. \textbf{133} (2005), 257--261.

\bibitem{Arv89}
W.~Arveson, \emph{{Continuous analogues of Fock space}}, Mem. Amer. Math. Soc.,
  no. 409, American Mathematical Society, 1989.

\bibitem{Bha96}
B.V.R. Bhat, \emph{{An index theory for quantum dynamical semigroups}}, Trans.\
  Amer.\ Math.\ Soc. \textbf{348} (1996), 561--583.

\bibitem{Ble97}
D.P. Blecher, \emph{{A new approach to Hilbert $C^*$--modules}}, Math.\ Ann.
  \textbf{307} (1997), 253--290.

\bibitem{GLR85}
P.~Ghez, R.~Lima, and J.~Roberts, \emph{{$W^*$--Categories}}, Pacific J. Math.
  \textbf{120} (1985), 79--109.

\bibitem{Kas80}
G.G. Kasparov, \emph{{Hilbert $C^*$--modules, theorems of Stinespring {\&}
  Voi\-cu\-les\-cu}}, J.\ Operator Theory \textbf{4} (1980), 133--150.

\bibitem{Lan95}
E.C. Lance, \emph{{Hilbert $C^*$--modules}}, Cambridge University Press, 1995.

\bibitem{MuSo00}
P.S. Muhly and B.~Solel, \emph{{On the Morita equivalence of tensor algebras}},
  Proc.\ London Math.\ Soc. \textbf{81} (2000), 113--168.

\bibitem{MuSo04}
\bysame, \emph{{Hardy algebras, $W^*$--correspondences and interpolation
  theory}}, Math.\ Ann. \textbf{330} (2004), 353--415.

\bibitem{MuSo05}
\bysame, \emph{{Duality of $W^*$-correspondences and applications}}, QP-PQ:
  Quantum Probability and White Noise Analysis XVIII (M.~Sch\"urmann and
  U.~Franz, eds.), World Scientific, 2005, pp.~396--414.

\bibitem{Pas73}
W.L. Paschke, \emph{{Inner product modules over $B^*$--algebras}}, Trans.\
  Amer.\ Math.\ Soc. \textbf{182} (1973), 443--468.

\bibitem{Rie74}
M.A. Rieffel, \emph{{Induced representations of $C^*$-algebras}}, Adv.\ Math.
  \textbf{13} (1974), 176--257.

\bibitem{Rie74a}
\bysame, \emph{{Morita equivalence for $C^*$--algebras and $W^*$--algebras}},
  J.\ Pure Appl.\ Algebra \textbf{5} (1974), 51--96.

\bibitem{Ske00b}
M.~Skeide, \emph{{Generalized matrix $C^*$--algebras and representations of
  Hilbert modules}}, Mathematical Proceedings of the Royal Irish Academy
  \textbf{100A} (2000), 11--38.

\bibitem{Ske01}
\bysame, \emph{{Hilbert modules and applications in quantum probability}},
  Ha\-bi\-li\-ta\-tions\-schrift, Cottbus, 2001, Available at \newline
  \tt{\footnotesize
  http://www.math.tu-cottbus.de/INSTITUT/lswas/\_skeide.html}.

\bibitem{Ske02}
\bysame, \emph{{Dilations, product systems and weak dilations}}, Math.\ Notes
  \textbf{71} (2002), 914--923.

\bibitem{Ske03c}
\bysame, \emph{{Commutants of von Neumann modules, representations of
  $\sB^a(E)$ and other topics related to product systems of Hilbert modules}},
  Advances in quantum dynamics (G.L. Price, B~.M. Baker, P.E.T. Jorgensen, and
  P.S. Muhly, eds.), Contemporary Mathematics, no. 335, American Mathematical
  Society, 2003, pp.~253--262.

\bibitem{Ske03p1}
\bysame, \emph{{Intertwiners, dual quasi orthonormal bases and
  representations}}, in preparation, 2004.

\bibitem{Ske04p}
\bysame, \emph{{~~Unit vectors, Morita equivalence and endomorphisms}}, Preprint,
  ArXiv: math.OA/0412231, 2004.

\bibitem{Ske05p}
\bysame, \emph{{Commutants of von Neumann Correspondences and Duality of
  Eilenberg-Watts Theorems by Rieffel and by Blecher}}, Preprint, ArXiv:
  math.OA/0502241, 2005, To appear in Proceedings of ``25th QP-Conference
  Quantum Probability and Applications'', Bedlewo 2004.

\bibitem{Ske05a}
\bysame, \emph{{Three ways to representations of $\sB^a(E)$}}, QP-PQ: Quantum
  Probability and White Noise Analysis XVIII (M.~Sch\"urmann and U.~Franz,
  eds.), World Scientific, 2005, pp.~504--517.

\bibitem{Tsi00p1}
B.~Tsirelson, \emph{{From random sets to continuous tensor products: answers to
  three questions of W.\ Arveson}}, Preprint, ArXiv: math.FA/0001070, 2000.

\end{thebibliography}

\newcommand{\Swap}[2]{#2#1}\newcommand{\Sort}[1]{}
\providecommand{\bysame}{\leavevmode\hbox to3em{\hrulefill}\thinspace}
\providecommand{\MR}{\relax\ifhmode\unskip\space\fi MR }
\providecommand{\MRhref}[2]{%
  \href{http://www.ams.org/mathscinet-getitem?mr=#1}{#2}
}
\providecommand{\href}[2]{#2}


\end{document}